\documentclass{article}
\usepackage{amsmath}
\usepackage{amssymb}
\usepackage{amsthm}
\newtheorem{theorem}{Theorem}
\newtheorem{corollary}{Corollary}
\newtheorem{lemma}{Lemma}
\newtheorem{Remark}{Remark}
\newtheorem{proposition}{Proposition}

\newcommand\dd{\mathsf{d}}
\newcommand\vkt[1]{{{\bf \bar{v}}}_{#1}^{\intercal}}
\newcommand\BB[2]{T_{#2}^{(#1)}}
\newcommand\baromega{\omega}
\newcommand\subo{A}
\newcommand\suboo{\underline{\omega}}
\newcommand\vkbt[1]{{{\bf u}}_{#1}^{\intercal}}
\newcommand\F{\mathsf{F}}

\newcommand\bomega{\bar{\omega}}

\newcommand\R{{\mathrm I}\! {\mathrm R}}
\def\f#1{{\textstyle \frac{1}{#1}}}
\def\fr#1#2{{\textstyle \frac{#1}{#2}}}
\newcommand\V{{\mathbf{v}}}
\def\sph#1{h^{(#1)}}
\def\vh{{\bf h}}
\newcommand\ku[2]{u_{#1}^{(#2)}}
\newcommand\kb[2]{{\bar{u}}_{#1}^{(#2)}}
\newcommand\vk[1]{{\bf u}_{#1}}
\newcommand\vkb[1]{{\bf \bar{u}}_{#1}}
\def\uu{\mathbf{u}}
\newcommand\Frame[1]{\mathbf{e}_{#1}}
\newcommand\bFrame[1]{\bar{\mathbf{e}}_{#1}}
\newcommand\ddt{\frac{\partial}{\partial \pr}}

\newcommand\ddx{\frac{\partial}{\partial \pt}}

\newcommand\MM{\mathsf{M}}
\newcommand\TT{\mathsf{T}}
\newcommand\torsion{\mathsf{T}}
\newcommand\cs[1]{\Frame{0 #1}}
\newcommand\E{{\mathcal E}}
\newcommand\cu{{\varkappa}}

\newcommand\pr{{t}}
\newcommand\pt{{x}}
\newcommand\lieg{\mathfrak{g}}
\newcommand\euc[1]{\mathfrak{euc}_{#1}}
\newcommand\orth[1]{\mathfrak{o}_{#1}}
\newcommand\lieh{\mathfrak{h}}

\newcommand\Hop{{\mathfrak{H}}}
\newcommand\D{{\mathcal D}}
\newcommand\Dx{D_\pt}
\newcommand\Dt{D_\pr}
\newcommand\Dxi{{D_\pt^{-1}}}
\newcommand\DX[1]{{D_\pt^{#1}}}

\newcommand\CC{{\mathcal C}}

\newcommand\N{{\mathbb N}}
\newcommand\Z{{\mathbb Z}}

\begin{document}
\bibliographystyle{alpha}
\title{Integrable Systems  in \(n\)-dimensional Riemannian Geometry}
\author{
Jan A. Sanders and Jing Ping Wang \\
Vrije Universiteit\\
Faculty of Sciences \\
Department of Mathematics\\
De Boelelaan 1081a\\
1081 HV Amsterdam\\
The Netherlands}
\date{\today}
\maketitle
\begin{abstract}
In this paper we show that if one writes down the structure
equations for the evolution of a curve embedded in
an \(n\)-dimensional Riemannian manifold with constant curvature
this leads to a symplectic, a Hamiltonian and an hereditary operator.
This gives us a natural connection between finite dimensional geometry,
infinite dimensional geometry and integrable systems.
Moreover one finds a Lax pair
in \(\orth{n+1}\) with the vector modified Korteweg-De Vries equation (vmKDV)
\[
\vk{t}=
 \vk{xxx}+\fr32 ||\vk{}||^2 \vk{x}
\]
as integrability condition.
We indicate that other integrable vector evolution equations
can be found by using a different Ansatz on the form of the Lax pair.
We obtain these results by using the {\em natural} or {\em parallel} frame
and we show how this can be gauged by a generalized
Hasimoto transformation to the (usual) {\em Fren{\^e}t} frame.
If one chooses the curvature to be zero, as is usual in the context
of integrable systems, then one loses information unless
one works in the natural frame.
\end{abstract}
\section{Introduction}
The study of the relationship between finite-dimensional
differential geometry  and partial differential equations
which later came to be known as integrable systems,
 started in the 19th century.
Liouville found and solved the equation describing
minimal surfaces in 3-dimensional Euclidean space \cite{Liou53}.
Bianchi solved the general Goursat problem for the sine-Gordon equation,
which arises in the theory of pseudospherical surfaces \cite{Bia10,Bia92}.

Much later Hasimoto \cite{Has72} found the relation between the equations for curvature and torsion
of vortex filament flow and the nonlinear Schr{\"o}dinger equation,
which led to many new developments,
cf.  \cite{MR85g:58039,MR92k:58118,MR95b:58075,MR96j:58082,MR99m:58119,MR1761235}.

The similarity between the structure equations for connections in differental geometry
and Lax pair equations for integrable systems has
intrigued many researchers from both fields
of mathematics. We refer the interested reader to the following books:
\cite{MR2002h:53003,MR99g:58036,MR2002f:53009,MR1913803,MR1908706}.
A good introductory review is \cite{MR98f:58111}.

Recently, we showed in \cite{MR1894465} that if
a flow of curves in a 3--dimensional Riemannian manifold with constant
curvature \(\cu\) follows an
arc--length preserving geometric evolution, the evolution of its
curvature and torsion
is always
a Hamiltonian flow with respect
to the pencil \(\E+\D+\cu \CC\),
where \(\E, \D\) and \(\CC\) are compatible Hamiltonian structures.
However, the close
geometric relationship remains here: the triplet is obtained
solely from the intrinsic geometry of curves on 3--dimensional
Riemannian manifolds with constant curvature.

Once one has two compatible Hamiltonian operators,
one can construct an hereditary operator, from which
a hierarchy of integrable equations can be computed.
They are all Hamiltonian with respect to two different
Hamiltonian operators, that is, biHamiltonian as defined in \cite{Mag78}.
Thus Poisson geometry is very
important in the study of integrable systems, cf.  \cite{MR94g:58260, MR94j:58081}.

The goal of the present paper is to generalize this analysis to
arbitrary dimension and see how much of the infinite dimensional geometric
structure is still present in this case.

The Hasimoto transformation is a Miura transformation,
which is induced by a gauge
transformation from the Fren{\^e}t frame (the obvious frame to choose
from the traditional differential geometry point of view) to the parallel
or natural frame, cf. \cite{MR95b:58075}.
On the basis of {\em loc.cit.}, Langer and Perline \cite{MR1781618} advocated
the use of the natural frame in the \(n\)-dimensional situation
in the context of vortex-filament flow.

We show in this paper that the appreciation
of the exact relationship between the underlying finite
dimensional geometry and the infinite dimensional geometry
has been complicated by the use
of the Fren{\^e}t frame.
In fact, the operators that lead to biHamiltonian systems
{\bf naturally} come out of the computation of the structure equations.
This fact is true in general, but only when one uses the parallel
or natural frame does this conclusion come out automatically.
This will be the first ingredient of our approach.

To relate our results to the classical situation in terms
of the curvatures of the evolving curve in the Fren{\^e}t frame we can of course
say that this relation is obvious since there must exist a gauge transform
between the two connections, defined by the Cartan matrices
specifying our frame.
Nevertheless, we give its explicit construction,
and produce a generalized Hasimoto transformation in arbitrary dimension,
with complete proof, as announced in \cite{Wang02}.
This insures that anything that can be formulated
abstractly can also be checked by direct (though complicated) computations.
Since the equations are derived from the computation of the curvature
tensor which behaves very nicely under transformations,
we feel that these direct calculations are the very last thing to try.
In this kind of problem the abstract point of view
seems to be much more effective computationally.
With hindsight one might say that Hasimoto was the first to exploit
the natural frame to show the integrability of the equations
for the curvature and torsion of a curve
embedded in a three-dimensional Riemannian manifold.

Here we should mention the fact that the generalization of the Hasimoto
also plays a role in \cite{MR1781618}, but it is in a different direction;
the focus is on the complex structure and this is
generalized to the \(2n\)-dimensional situation.

A second ingredient in our approach is the fact that we assume the
Riemannian manifold to have (nonzero) constant curvature.
If one takes the curvature equal to zero and one uses the natural frame,
then one can still recover the symplectic and cosymplectic operators
which give rise to a biHamiltonian system.
But if one uses another frame, then information gets lost
and one is faced with the difficult task of recreating this information
in order to get the necessary operators and Poisson geometry.
In our opinion, this has been a major stumbling block in the analysis
of the relation between high-dimensional geometry and integrability.

A third ingredient is the study of symmetric spaces.
The relation between symmetric spaces and integrability was
made explicit in the work of Fordy and others \cite{MR88j:58046,MR84k:58106},
There it is shown how to construct a Lax pair using the structure
of the Lie algebra. This construction does not work in the semidirect
product \(\euc{n}\), but
there we use the fact that we can map the elements in the
\(n\)-dimensional euclidean algebra to the \(n+1\)-dimensional
orthogonal algebra.
This leads immediately to the construction of Lax pairs, combining
the methods in \cite{MR2002k:37141} and \cite{MR1781618}.

Putting together these three elements, natural frame, constant curvature
and extension to a symmetric semisimple Lie algebra,
gives us complete control over all the integrability issues one might like
to raise. We have tried to formulate things in such a way that it is clear
what to do in other geometries. The theory, however, is not yet strong
enough to guarantee success in other geometries, so this will
remain an area of future research.

If we try to understand the success of the natural frame in this context,
then it seems that the key is the natural identification
of the Lie algebra \(\euc{n}/\orth{n}\)
with the adjoint orbit of \(\orth{n}^{*}\).
In the natural frame, this identification is simple to see
and can be performed by putting a complex structure on \(\orth{n+1}\).
To give a geometrical explanation of this, one may have to
formulate everything in terms of Poisson reduction of the Kac-Moody
bracket of \(SO(n)\), but we make no attempt
to do so here.
The main point we will make in this paper is that the whole
construction can be explicitly computed. Once this is done, one
can then choose a more conceptional approach later on. We have
indicated the general character of the approach where appropriate.

The paper is organized as follows.
In section \ref{Sec1} we derive the structure equations,
using Cartan's moving frame method,
and find that the equation is
of the type
\[
\vk{t}=\Phi_1 \vh +\cu \Phi_2 \vh, \quad \Phi_2=Id_{n-1}.
\]
In section \ref{Sec3} we show that \(\Phi_1\) is an hereditary operator,
and we draw the conclusion that if we had been computing in any other
frame, resulting in the equation
\begin{eqnarray}\label{flow}
\vkb{t}=\bar{\Phi}_1 \bar{\vh}+\cu \bar{\Phi}_2 \bar{\vh},
\end{eqnarray}
then \(\bar{\Phi}_1 \bar{\Phi}_2^{-1}\) is also hereditary.

We then proceed, by taking \(\vh=\vk{1}\), the first \(x\)-derivative of \(\vk{}\),
to derive a vector mKDV equation.
In \cite{MR95b:58075} it was predicted (or derived by general considerations)
that the equation should be a vmKDV equation,
but there are two versions of vmKDV equation, cf. \cite{MR1872986},
and the prediction does not say which version it should be.
That question is now settled.

In section \ref{secHasimoto} we state explicitly the formula for
the generalized Hasimoto transformation, which transforms the
Fren{\^e}t frame into the natural frame and in Appendix
\ref{Hasimoto} we prove this.

We conclude in section \ref{Laxpair} by constructing a Lax pair which has this equation
as its integrability condition.
This is done as follows.
The geometric problem is characterized by two Lie algebras,
\(\orth{n}\) and  \(\euc{n}\).
The usual Lie algebraic construction of Lax pairs starts with a
symmetric Lie algebra \(\lieh^0+\lieh^1\),
with
\begin{eqnarray}\label{sym}
[\lieh^0,\lieh^0]\subset \lieh^0, [\lieh^0,\lieh^1]\subset \lieh^1,
[\lieh^1,\lieh^1]\subset \lieh^0.
\end{eqnarray}
If we think of \(\orth{n}\) as \(\lieh^0\),
then we cannot take \(\euc{n}\) as \(\lieh^0+\lieh^1\),
since \(\lieh^1\) is then abelian. Although it trivially obeys all the inclusions,
this triviality kills all the actions we need to construct the Lax pair.
But if we take \(\orth{n+1}
=\lieh^0+\lieh^1\),
then things fit together nicely and we can simply copy the construction
in \cite{MR2002k:37141}.
Of course, we can identify \(\lieh^1\) with \(\euc{n}/\orth{n}\) as vectorspaces,
and this should be the main consideration (apart from the requirements
in (\ref{sym})) in other geometries when
one tries to find such a {\em symmetric extension} of the geometrically
given \(\lieh^0\).

{\bf Acknowledgment}
The authors thank Gloria Mar\'\i\ Beffa for numerous discussions,
which both influenced the questions that we asked ourselves as well
as some the methods we used.
The work was mainly carried out at the Isaac Newton Institute
for Mathematical Sciences
during the Integrable Systems programme in 2001 and
we thank $\mathcal{NWO}$ (Netherlands Organization for Scientific Research)
for financial support and the Newton Institute for its
hospitability and support.
Thanks go to Hermann Flaschka for reminding us of \cite{MR95b:58075}
and to Sasha Mikhailov for his help with Lax pairs.
\section{Moving Frame Method}\label{Sec1}
Our approach, although restricted to the Riemannian case,
is also suitable for arbitrary geometries. As a general reference
to the formulation of geometry in terms of Lie algebras, Lie groups and connections
we refer to \cite{MR98m:53033}.

We consider a curve (denoted by \(\Frame{0}\) if we want to stress its
relation to the other vectors in the moving frame,
or \(\gamma\), if we just want to concentrate on the curve itself),
parametrized by arclength \(x\) and evolving
geometrically in time \(t\),
which is imbedded in a Riemannian manifold \(\MM\).
We choose a frame \(\Frame{1},\cdots,\Frame{n}\)
as a basis for \(\TT_{\Frame{0}}\MM\),
and a dual basis of one forms \(\tau_1,\cdots,\tau_n\in\bigwedge^1 \MM\),
so that \(\tau_i(\Frame{j})=\delta_{i}^{j}\).
Choosing a frame is equivalent to choosing an element in \(g\in G=Euc(n,\R)=O(n,\R)\ltimes \R^n\).
We write \(g=\left(\begin{array}{cc}1 & \Frame{0}\\0&\Frame{}
\end{array}\right)\).
Let  \(\Gamma(\TT\MM)\) be the space of smooth sections,
that is, maps \(\sigma: \MM\rightarrow \TT\MM\) such that
\(\pi \sigma\) is the identity map on \(\MM\),
where \(\pi \) is the natual projection of the tangent space
on its base \(\MM\).
We consider the \(\Frame{i}\) as sections in \(\TT\MM\), at least locally,
that is, elements in \(\Gamma(\TT\MM)\), varying \(\Frame{0}\).
We assume that there is a Lie algebra
\(\lieg\) and a subalgebra \(\lieh\)
such that \(\TT_{\Frame{0}} \MM \backsimeq \lieg/\lieh\) as vectorspaces.

We now define a connection
\(\dd: \bigwedge^p\MM\otimes\Gamma(\TT\MM)\rightarrow \bigwedge^{p+1}
\MM \otimes \Gamma(\TT\MM)\)
as follows.
Let \(\omega\in \bigwedge^p\MM\). Then
\[
\dd \omega\otimes\sigma=d\omega\otimes \sigma+{(-1)}^p\omega\otimes\dd\sigma.
\]
Notice that we use
\( d\) as a cochain map from \(\bigwedge^p\MM\)
to \(\bigwedge^{p+1}\MM\), \(p\in\N\) in the ordinary de Rham complex.
We extend the connection to act on \(\Frame{0}\) as follows.
\[\dd\Frame{0}=\sum_i \tau_i\otimes \Frame{i},\quad \tau_i\in\bigwedge^1\MM
\]
with \(\tau=(\tau_1,\cdots,\tau_n)\in C^1(\lieg/\lieh,\lieg/\lieh)\)
a differential \(1\)--form.
We now suppose that \(\dd \Frame{i}=\sum_{j} \omega_{ij} \otimes\Frame{j}\),
with \(\omega_{ij}\in\bigwedge^1\MM\), that is,
\(\omega\in C^1(\lieg/\lieh,\lieh)\).
This connection naturally extends to a connection
\(\dd: \bigwedge^p\MM\otimes\Gamma(G)\rightarrow \bigwedge^{p+1}
\MM \otimes \Gamma(G)\).
We may write our basic structure equations as
\begin{eqnarray}\label{MC}
\dd g=A g
=
\left(\begin{array}{cccc}0&\tau_{}
\\0& \omega
\end{array}\right)
g
\end{eqnarray}
where \(A\) is the vector potential or Cartan matrix.
 We remark
that \(d\tau(X,Y)=X\tau(Y)-Y\tau(X)-\tau([X,Y])\). Differentiating
both sides of (\ref{MC}), we have
\begin{eqnarray*}
\dd^2
g
&=&
\left(\begin{array}{cc}0&d\tau_{}
\\0& d\omega
\end{array}\right)
g
-
\left(\begin{array}{cc}0&\tau_{}
\\0& \omega
\end{array}\right)
\left(\begin{array}{cc}0&\dd\Frame{0}\\0&\dd\Frame{}
\end{array}\right)
\\&=&
\left(\begin{array}{cc}0&d\tau_{}
\\0& d\omega
\end{array}\right)
g
-
\left(\begin{array}{cc}0&\tau_{}
\\0& \omega
\end{array}\right)
\left(\begin{array}{cc}0&\tau_{}
\\0& \omega
\end{array}\right)
g
\\&=&
\left(\begin{array}{cc}0&d\tau_{}-\tau\wedge\omega
\\0& d\omega-\omega\wedge\omega
\end{array}\right)
g
=\F g
\end{eqnarray*}
We  draw the conclusion that
\[
\dd^2=\F=
\left(\begin{array}{cc}0&\torsion
\\0& \Omega
\end{array}\right)
=\left(\begin{array}{cccc}0&d\tau_{}-\tau\wedge\omega
\\0&d\omega-\omega\wedge\omega
\end{array}\right),
\]
where \(\F\in C^2(\lieg/\lieh,\lieg)\), of which
the \(2\)--form \(\Omega\in C^2(\lieg/\lieh,\lieh)\) is called the {\em curvature},
the \(2\)--form \(\torsion\) is called {\em torsion}
and \(\omega\wedge \omega\) is often denoted by \(\frac{1}{2} [\omega,\omega]\).

We are now in a position to compute a few things explicitly.
Let \(X=\cs{\ast}\ddx=\Dx\) and \(Y=\cs{\ast}\ddt=\Dt\),
that is, the induced vectorfields tangent to the imbedded curve
of \(\ddx\) and \(\ddt\), respectively.

Then, using the fact that \([X,Y]=0\), we obtain
\begin{eqnarray}
\lefteqn{\F(X,Y)=}&&\nonumber
\\&=&\left(\begin{array}{cc}0&\torsion(X,Y)
\\0& \Omega(X,Y)
\end{array}\right)
=\left(\begin{array}{cccc}0&d\tau_{}(X,Y)-\tau(X)\wedge\omega(Y)
\\0& d\omega(X,Y)-\omega(X)\wedge\omega(Y)
\end{array}\right)\nonumber
\\&=&
\left(\begin{array}{cccc}0&X\tau(Y)-Y\tau(X)-\tau(X)\wedge\omega(Y)
\\0&
X\omega(Y)-Y\omega(X)-\omega(X)\wedge\omega(Y)
\end{array}\right)\nonumber
\\&=&
\left(\begin{array}{cccc}0&\Dx\tau(Y)-\Dt\tau(X)-\tau(X)\wedge\omega(Y)
\\0&
\Dx\omega(Y)-\Dt\omega(X)-\omega(X)\wedge\omega(Y)
\end{array}\right)\label{CuE}
\\&=&
[\Dx+\left(\begin{array}{cccc}0&\tau(X)
\\0&
\omega(X)
\end{array}\right),
\Dt+\left(\begin{array}{cccc}0&\tau(Y)
\\0& \omega(Y) \end{array}\right)]\label{zero}
\end{eqnarray}
\section{Hasimoto transformation}\label{secHasimoto}
\renewcommand\vk[1]{{\bf \bar{u}}_{#1}}
\renewcommand\vkb[1]{{\bf u}_{#1}}
\renewcommand\kb[2]{u_{#1}^{(#2)}}
\renewcommand\ku[2]{{\bar{u}}_{#1}^{(#2)}}
\renewcommand\bFrame[1]{\mathbf{e}_{#1}}
\renewcommand\Frame[1]{\bar{\mathbf{e}}_{#1}}

In 1975, Bishop discovered the same transformation as Hasimoto
when he studied the relations between two different frames to
frame a curve in 3-dimensional Euclidean space, namely the
Fren{\^e}t frame and parallel (or natural) frame, cf.
\cite{MR51:6604}. More explicitly, let the orthonormal basis
\(\{T, N, B\}\) along the curve be the Fren{\^e}t frame, that is,
\begin{eqnarray*}
\left( \begin{array}{l} T\\ N\\ B
\end{array} \right)_x=\left( \begin{array}{lll}
\ 0 & \ \kappa & 0\\ -\kappa & \ 0 & \tau\\ \ 0 & -\tau & 0
\end{array} \right)\left( \begin{array}{l}
T\\ N\\ B \end{array} \right).
\end{eqnarray*}
The matrix in this equation is called the Cartan matrix. Now we
introduce the following new basis
\begin{eqnarray}\label{Tf}
\left( \begin{array}{l} T\\N^1\\ N^2
\end{array} \right)=\left( \begin{array}{lll}
1\ & 0& \, 0\\0\ & \cos\theta & -\sin \theta \\ 0\ & \sin \theta &
\, \cos \theta
\end{array} \right)\left( \begin{array}{l}
T\\N\\ B \end{array} \right),\quad \theta=\int \tau dx\ .
\end{eqnarray}
Its frame equation is
\begin{eqnarray*}
\left( \begin{array}{l} T\\ N^1\\ N^2
\end{array} \right)_x=\left( \begin{array}{lll}
0 & \kappa \cos \theta& \kappa \sin \theta\\
-\kappa \cos \theta & \ \ 0 & \ \ 0 \\ -\kappa \sin \theta & \ \ 0
& \ \ 0
\end{array} \right)\left( \begin{array}{l}
T\\ N^1\\ N^2 \end{array} \right).
\end{eqnarray*}
We call the basis \(\{T, N^1, N^2\}\) the parallel frame. This
geometric meaning of Hasimoto transformation was also pointed out
in \cite{MR95b:58075}, where the authors studied the connection
between the differential geometry and integrability.

In this section, we give explicit formula of such a transformation
in arbitrary dimension \(n\), which has the exact same geometric
meaning as the Hasimoto transformation. Therefore, we call it
generalized Hasimoto transformation. The existence of such a
transformation is clear from the geometric point of view, and was
mentioned and implicitly used in several papers such as
\cite{MR51:6604} and \cite{MR1781618}.

First we give some notation. Denote the Cartan matrix of the
Fren{\^e}t frame  \(\Frame{}\) by \(\bomega(\Dx)\) (
\(\Frame{x}=\bomega(\Dx) \Frame{}\) ) and that of the parallel
frame \(\bFrame{}\) by \(\baromega(\Dx)\) (
\(\bFrame{x}=\baromega(\Dx) \bFrame{}\) ), i.e.,
\begin{eqnarray*}
\bomega(\Dx)=\left( \begin{array}{ccccc} 0&\ku{}{1}&0&\cdots&0
\\
-\ku{}{1}&0&\ku{}{2}&\cdots&0
\\
\vdots&\ddots&\ddots&\ddots&\vdots
\\
0&-\ku{}{n-3}&0&\ku{}{n-2}&0
\\
0&0&-\ku{}{n-2}&0&\ku{}{n-1}
\\
0&0&0&-\ku{}{n-1}&0
\end{array}\right)
\end{eqnarray*}
and, letting \(\vkb{}=(\kb{}{1},\ \kb{}{2},\cdots  \kb{}{n-1})^\intercal\),
\begin{eqnarray*}
\baromega(\Dx)=\left( \begin{array}{ll} 0 & \vkbt{}\\
-\vkb{} & 0
\end{array} \right).
\end{eqnarray*}
The orthogonal matrix that keeps the first row and rotates
the \(i\)-th and \(j\)-th row with the angle \(\theta_{ij}\)
is denoted by \(R_{ij}\), where \(2\leq i<j\leq n\). For example,
when \(n=3\), the orthogonal matrix
\begin{eqnarray*}
R_{23}= \left(\begin{array}{lll} 1 & \ 0 & 0\\
0&\  \cos \theta_{23} & \sin \theta_{23}\\
0& -\sin \theta_{23} & \cos \theta_{23}\end{array}\right)
\end{eqnarray*}
is the transformation between two frames, i.e., \(\Frame{}= R_{23}
\bFrame{}\), comparing to (\ref{Tf}). Let us compute the
conditions such that
\begin{eqnarray*}
\bomega(\Dx) R_{23}- \frac{\partial R_{23}}{\partial x}= R_{23}
\baromega(\Dx).
\end{eqnarray*}
Then we have \(\ku{}{2}=\Dx \theta_{23}\), \(\kb{}{1}=\ku{}{1}\cos \theta_{23}\)
and \(\kb{}{2}=\ku{}{1}\sin \theta_{23}\).
This is the famous Hasimoto transformation, i.e.,
\(\phi=\kb{}{1}+i \kb{}{2}=\ku{}{1} \exp ( i \int \ku{}{2} dx)\).

Therefore, the generalized Hasimoto transformation,
which by definition transforms the Fren{\^e}t frame
into the natural frame, can be found by
computing the orthogonal matrix, which
gauges the Cartan matrix of the Fren{\^e}t frame into
that of the parallel frame in \(n\)-dimensional space.

\begin{theorem}\label{TheoHasimoto}
Let \(n\geq 3\). The orthogonal matrix \(R\)
gauges the standard Fren{\^e}t frame into the natural frame,
that is,
\[
\bomega(\Dx) R-\Dx R =R \baromega(\Dx),
\]
where \(R=R_{n-1,n}\cdots R_{3,n} \cdots R_{34}R_{2,n}\cdots
R_{24} R_{23}\). This leads to the components of \(\uu\)
satisfying the Euler transformation
\begin{eqnarray*}
\left\{ \begin{array}{lll}
 \kb{}{1}=
\ku{}{1} \cos \theta_{23} \cos \theta_{24}\cdots \cos \theta_{2n},\\
\kb{}{2}=
\ku{}{1} \sin \theta_{23} \cos \theta_{24}\cdots \cos \theta_{2n}, \\
\kb{}{3}=
\ku{}{1} \sin \theta_{24}\cdots \cos \theta_{2n}, \\
\cdots \cdots \\
\cdots \cdots \\
 \kb{}{n-2}=
\ku{}{1} \sin \theta_{2,n-1} \cos \theta_{2n},\\
 \kb{}{n-1}=
\ku{}{1} \sin \theta_{2n},
\end{array} \right.
\end{eqnarray*}
and the curvatures in standard Fren{\^e}t frame satisfying
\begin{eqnarray}\label{AFormu}
&&\ku{}{i}=\prod_{j=i+2}^n \frac{\cos \theta_{i,j}}{\cos \theta_{i+1,j}} \subo_i,
\quad i=2,3,\cdots, n-1,
\end{eqnarray}
where the \(\subo_i\) are generated by the recursive relation
\begin{eqnarray}\label{AFormA}
\subo_i=\Dx \theta_{i,i+1}+\sin \theta_{i-1,i+1} \subo_{i-1}, \ \subo_1=0, \ 2\leq i\leq n-1.
\end{eqnarray}
There are \(\frac{1}{2} (n-2)(n-3)\) constraints on the rotating angles:
\begin{eqnarray}\label{AFormth}
&&\Dx \theta_{ij}=\prod_{l=i+2}^j\frac{\cos \theta_{il}}{\cos \theta_{i+1,l}}
\sin \theta_{i+1,j} \subo_i
-\prod_{l=i+1}^{j-1}\frac{\cos \theta_{i-1,l}}{\cos \theta_{i,l}}
\sin \theta_{i-1,j} \subo_{i-1},
\end{eqnarray}
where \(4\leq j\leq n\), \(2\leq i\leq j-2\).
\end{theorem}
\begin{proof}
It is obvious that \(R= \left( \begin{array}{ll} 1 & 0\\ 0 & T
\end{array} \right)\), where \(T\) is an orthogonal
\((n-1)\times(n-1)\)-matrix and the first row of \(T\) equals is
\(\frac{\vkbt{}}{\ku{}{1}}\), that is,
\[\left(\cos \theta_{23} \cos \theta_{24}\cdots \cos \theta_{2n},
\sin \theta_{23} \cos \theta_{24}\cdots \cos \theta_{2n}, \cdots
,\sin \theta_{2n}\right).\] This implies that \(\vkt{}
T=\vkbt{}\), where \(\vkt{}=(\ku{}{1},0, \cdots, 0)\).
\renewcommand\vk[1]{\bar{\bf v}_{#1}}
Since \(T T^\intercal=I\), we have \(T  \vkb{}=\vk{}\). Rewrite
\(\bomega(\Dx)=\left( \begin{array}{ll} 0 & \vkt{}\\ -\vk{} &
\suboo
\end{array} \right)\), and compute
\begin{eqnarray*}
&&\bomega(\Dx) R-\Dx R - R \baromega(\Dx) \\
&=&\left( \begin{array}{ll} 0 & \vkt{}\\ -\vk{} & \suboo
\end{array} \right) \left( \begin{array}{ll} 1 & 0\\ 0 & T
\end{array} \right) -\left( \begin{array}{ll} 0 & 0\\ 0 & \Dx T
\end{array} \right) -\left( \begin{array}{ll} 1 & 0\\ 0 & T
\end{array} \right) \left( \begin{array}{ll} 0 & \vkbt{}\\ -\vkb{}
& 0 \end{array} \right)
\\
&=&\left( \begin{array}{ll} 0 & \vkbt{}\\ -\vk{} &\suboo T
\end{array} \right) -\left( \begin{array}{ll} 0 & 0\\ 0 & \Dx T
\end{array} \right) -\left( \begin{array}{ll} 0 & \vkbt{}\\ -\vk{}
& 0 \end{array} \right)
\\
&=&\left( \begin{array}{ll} 0 & 0\\ 0 &\suboo T-\Dx T \end{array}
\right)
\end{eqnarray*}
Therefore, to prove the statement, we only need to solve \(\suboo
T=\Dx T\), that is, the matrix \(T\) gauges \(\suboo\) into zero,
which is proved in Appendix \ref{Hasimoto}.
\end{proof}

\renewcommand\ku[2]{{u}_{#1}^{(#2)}}
\renewcommand\vk[1]{{\bf u}_{#1}}
\renewcommand\vkb[1]{{\bf \bar{u}}_{#1}}
\renewcommand\kb[2]{\bar{u}_{#1}^{(#2)}}

\renewcommand\Frame[1]{\mathbf{e}_{#1}}
\renewcommand\bFrame[1]{\bar{\mathbf{e}}_{#1}}
\section{Hereditary operator in Riemannian Geometry}\label{Sec3}
In this section, we present a hereditary operator that arises in
a natural way from the geometric arelength-preserving evolution of
curves in \(n\)-dimensional Riemannian manifold with constant curvature.
The operator is the product of a cosymplectic (Hamiltonian) operator
and a symplectic operator.

\begin{theorem}\label{Th1}
Let \(\gamma(\pt,\pr)\) be a family of curves on \(\MM\)
satisfying a geometric evolution system of equations of the form
\begin{equation}\label{uev}
\gamma_{\pr}= \sum_{l=1}^n \sph{l} \Frame{l}
\end{equation}
where \(\{ \Frame{l}, l=1,\cdots, n\}\) is the natural frame
of \(\gamma\), and where \(h_l\)
are arbitrary smooth functions of the curvatures
\(\ku{}{i}, i=1,\cdots n-1\) and their derivatives
with respect to \(x\).

Assume that \(\pt\) is the arc-length parameter and that
evolution (\ref{uev}) is arc-length preserving.

Then, the curvatures \(\vk{}=(\ku{}{1},\cdots,
\ku{}{n-1})^\intercal\) satisfy the evolution
\begin{eqnarray*}
  \vk{t}=\Re \vh{}-\cu \vh{},\quad \vh{}=(\sph{2},\cdots,\sph{n})^\intercal,
\end{eqnarray*}
where the operator \(\Re\) is hereditary and defined as follows:
\begin{eqnarray*}
\Re&=&\Dx^2+\langle \vk{},\ \vk{}\rangle +\vk{1} \Dxi \langle \vk{},\ \cdot\rangle
-\sum_{i<j}^{n} J_{ij} \vk{} \Dxi
\langle J_{ij} \vk{1},\ \cdot\rangle ,
\end{eqnarray*}
and where the \(J_{ij}\) are anti-symmetric matrices with nonzero entry of \((i,j)\)
being \(1\) if \(i<j\), that is, \((J_{ij})_{kl}=\delta_k^i\delta_j^l-\delta_l^i\delta_k^j\).
Moreover, \(\Re\) can be written as \(\mathfrak{HI}\), where
\(\mathfrak{I}\) is a symplectic operator defined by
\(\mathfrak{I}=\Dx+\uu \Dxi\uu^\intercal\),
and \(\mathfrak{H}\) is a cosymplectic (or Hamiltonian) operator, defined by
\begin{eqnarray*}
\mathfrak{H}=\Dx+\sum_{i=1}^{n-1}\sum_{j=i+1}^n
J_{ij}\uu\Dxi(J_{ij}\uu)^\intercal.
\end{eqnarray*}
\end{theorem}
\begin{proof}
In the Riemannian case, \(\lieh\) is \(\orth{n}\) and \(\lieg\) is \(\euc{n}\).
By fixing the frame, we know the value of \(\omega(\Dx)\) and \(\tau(\Dx)\).
We assumed our frame to be the natural (or parallel) frame,
see \cite{MR51:6604}.
So we have \(\Frame{1}=X\).  Due to \(\dd \Frame{0}=\tau \otimes\Frame{}\), we know
\[\tau(\Dx)=(1,0,\cdots, 0).\]
The Cartan matrix of the natural frame is by definition
\[
\omega(\Dx)=\left(\begin{array}{ccccc} 0 & \ku{}{1} &
\ku{}{2}&\cdots&
\ku{}{n-1}\\
-\ku{}{1} &0&0&\cdots&0\\
\vdots&\vdots&\vdots&\vdots&\vdots\\
-\ku{}{n-1}&0&0&\cdots&0
\end{array}\right),
\]
that is,
\(
\omega_{ij}(\Dx)=\delta_1^i \ku{}{j-1}-\delta_1^j
\ku{}{i-1}.
\)

Now we use the moving frame method in section \ref{Sec1} to derive
the evolution equation of the curvatures of a smooth curve \(\gamma(\pt,\pr)\)
in \(n\)-dimensional Riemannian manifold satisfying a geometric evolution of the form
\begin{eqnarray}\label{Curve}
\gamma_{\pr}= \sum_{i=1}^n \sph{i} \Frame{i},
\end{eqnarray}
where \(\{ \Frame{i}, i=1\cdots n\}\) is the natural frame and
\(\sph{i}\) are arbitrary smooth functions of the curvatures
\(\ku{}{i}, i=1\cdots n-1\), and their derivatives with respect to
the arclength parameter \(x\).

The given curve (\ref{Curve}) implies \(\sph{i}=\tau_i(\Dt)\).
So the equation (\ref{CuE}) reads
\begin{eqnarray*}
\torsion_{i}(\Dx,\Dt)&=&
\Dx \sph{i}-\omega_{1i}(\Dt)+ \ku{}{i-1}\sph{1}
-\delta_1^i\sum_{j=2}^n \sph{j} \ku{}{j-1}
\end{eqnarray*}
and we find
\begin{eqnarray}\label{Tor}
\omega_{1i}(\Dt)=\Dx \sph{i}+ \ku{}{i-1}\sph{1}
-\delta_1^i\sum_{j=2}^n \sph{j} \ku{}{j-1}
-\sum_{j=2}^n \torsion_{i}(\Dx,\Frame{j})\sph{j}.
\end{eqnarray}
We redefine \(\omega\) by
\begin{eqnarray*}
\tilde{\omega}_{1i}=\omega_{1i}+\torsion_{i}(\Dx,\cdot).
\end{eqnarray*}
This does not influence our previous calculations
since \(\omega(\Dx)=\tilde{\omega}(\Dx)\), and
reflects the fact that we can define the torsionless connection,
i.e., Riemannian connection. It is customary to call the geometry Riemannian
if the torsion is zero.
We now write \(\omega_{1i}\)
for \(\tilde{\omega}_{1i}\) to avoid the complication of the notation
and rewrite (\ref{Tor}) as
\begin{eqnarray}\label{Torless}
\omega_{1i}(\Dt)&=& \Dx \sph{i}+ \ku{}{i-1}\sph{1}
-\delta_1^i\sum_{j=2}^n \sph{j} \ku{}{j-1}.
\end{eqnarray}
We need \(\omega_{11}=0\) for the consistence, i.e.,
\begin{eqnarray*}
0&=& \Dx \sph{1} -\sum_{j=2}^n \sph{j} \ku{}{j-1}.
\end{eqnarray*}
Geometrically it means that the evolution is arc-length preserving.
By eliminating \(\sph{1}\) in (\ref{Torless}), we obtain,
taking the integration constants eequal to zero,
\[
\omega_{1i}(\Dt)=\Dx \sph{i} +\ku{}{i-1}\Dxi \sum_{j=2}^n \sph{j}
\ku{}{j-1}
=\mathfrak{I} \vh,\quad 1<i\leq n,
\]
with \(\vh=(h_2,\cdots,h_n)\).
This defines a symplectic operator \(\mathfrak{I}=\Dx+\uu \Dxi\uu^\intercal\),
cf. \cite{Wang02}, as will be proved in Proposition \ref{SymplecticOp} in
Appendix \ref{OpProof}.
We now do the same for \(\orth{n}\)-component in (\ref{CuE}).
Assuming \(j>i\) (\(\omega\) is antisymmetric), we have
\begin{eqnarray*}
\lefteqn{\Omega_{ij}(\Dx,\Dt)=
\Dx\omega_{ij}(\Dt)-\Dt (\delta_1^i \ku{}{j-1}-\delta_1^j
\ku{}{i-1})
}\\&&
-\sum_{l=1}^n (\delta_1^i \ku{}{l-1}-\delta_1^l
\ku{}{i-1}) \omega_{lj}(\Dt)
+\sum_{l=1}^n \omega_{il}(\Dt) (\delta_1^l \ku{}{j-1}-\delta_1^j
\ku{}{l-1}).
\end{eqnarray*}
This leads to
\begin{eqnarray*}
\Dt  \ku{}{j-1} &=& \Dx\omega_{1j}(\Dt) -\sum_{l=2}^n
\ku{}{l-1} \omega_{lj}(\Dt)
-\Omega_{1j}(\Dx,\Dt),\quad j=2,\cdots,n
\end{eqnarray*}
when \(i=1\), and when \(i>1\),
\begin{eqnarray*}
\Dx \omega_{ij}(\Dt) =
\Omega_{ij}(\Dx,\Dt)
- \ku{}{i-1} \omega_{1j}(\Dt)
+ \ku{}{j-1} \omega_{1i}(\Dt) .
\end{eqnarray*}
We combine these to
\begin{eqnarray*}
\lefteqn{\Dt  \ku{}{j-1} = -\Omega_{1j}(\Dx,\Dt)+\Dx\omega_{1j}(\Dt) }&&
\\&&
-\sum_{l=2}^n \ku{}{l-1} \Dxi (-
\ku{}{l-1} \omega_{1j}(\Dt)
+ \ku{}{j-1} \omega_{1l}(\Dt)
)
\quad j=2,\cdots,n
\end{eqnarray*}
Assume that the curvature of Riemannian manifold is constant,
i.e., the only nonzero entries in \(\Omega(\Frame{i}, \Frame{j})\) are
\[\Omega_{ij}(\Frame{i}, \Frame{j})=-\Omega_{ji}(\Frame{i}, \Frame{j})=\cu,
\quad j>i.\]
Then \(\Omega_{1j}(\Dx,\Dt)=h_j\).
The above formula defines in the case of constant curvature a cosymplectic
(or Hamiltonian) operator \(\mathfrak{H}\) (cf. \cite{Wang02}) given by
\begin{eqnarray*}
\mathfrak{H}=\Dx+\sum_{i=1}^{n-1}\sum_{j=i+1}^n
J_{ij}\uu\Dxi(J_{ij}\uu)^\intercal,
\end{eqnarray*}
where the \(J_{ij}\), given by
\((J_{ij})_{kl}=\delta_{ik}\delta_{jl}-\delta_{il}\delta_{jk}\),
are the generators of \(\orth{n}\)
and such that \(\Dt\uu=\mathfrak{H}\omega_{1\cdot}(\Dt)\).
The fact that the operator is indeed Hamiltonian will be given in Proposition
\ref{HamiltonianOp} in Appendix \ref{OpProof}.
The operator is weakly non-local, as defined in \cite{MR2002g:37093}.

It corresponds to what is called the {\em second Poisson operator \(P_u\)},
defined by
\[
P_u=\Dx+\pi_1 ad(u)+ad(u)\Dxi \pi_0 ad(u),
\]
in \cite{MR2002k:37141}, but looks simpler by the choice of frame.

We write the evolution equation for the curvatures \(\ku{}{j}\) as
\begin{eqnarray*}
\vk{t}=\Re \bf{h}-\cu \bf{h},
\end{eqnarray*}
where the operator \(\Re=\mathfrak{H}\mathfrak{I}\) is explicitly given as
\begin{eqnarray}\label{Rop}
\Re&=&\Dx^2+\langle \vk{},\ \vk{}\rangle+\vk{1} \Dxi \langle \vk{},\ \cdot\rangle
-\sum_{i<j}^{n} J_{ij} \vk{} \Dxi
\langle J_{ij} \vk{1},\ \cdot\rangle.
\end{eqnarray}
We prove that \(\Re\) is hereditary in Proposition \ref{HereditaryOp} in Appendix \ref{OpProof}.
\end{proof}
\begin{Remark}
The recursion operator can be used to generate infinitely many compatible
symplectic and cosymplectic weakly non-local operators starting
with \(\mathfrak{I}\) and \(\Hop\). That the result is indeed
weakly non-local follows from the techniques in \cite{SW00b,MR2002g:37093}.
\end{Remark}
\begin{Remark}
A special case arises for \(n=2\): one has
\[\Re=-(J_{12}(\Dx+\uu\Dxi\uu))^2=-R_{NLS}^2,\]
where \(R_{NLS}\) is the recursion
operator of the Nonlinear Schr{\"o}dinger (NLS) equation.
\end{Remark}
Let us now take \(\sph{j}= \ku{1}{j-1}\) in Theorem \ref{Th1}. Then
we obtain
\begin{eqnarray*}
 \vk{t} =\vk{3} +\frac{3}{2}
\langle  \vk{},\ \vk{}\rangle  \vk{1}
-\cu \vk{1}\mbox{   (vmKDV)},
\end{eqnarray*}
which is one of two versions of vector mKDV equations appearing in
\cite{MR1872986}.
\begin{corollary}
The hereditary operator defined in (\ref{Rop}) is a recursion operator of the vector
mKDV equation
\(
 \vk{t} =\vk{3} +\frac{3}{2} \langle  \vk{},\ \vk{}\rangle  \vk{1}.
\)
\end{corollary}
Although we give the explicit formula of the generalized Hasimoto
transformation in section \ref{secHasimoto}, it is difficult to
use it to compute the formula (\ref{flow}). It is easier to obtain
it directly using the moving frame method in section \ref{Sec1}.
\begin{theorem}
Let us assume that in the Fren{\^e}t frame, the curvatures of a
geometric curve \( \gamma_{\pr}= \sum_{l=1}^n \bar{\sph{l}}\bFrame{l} \) 
under arc-length preserving satisfy
\begin{eqnarray}\label{flow1}
\vkb{t}=\bar{\Phi_1}\bar{\vh}-\cu\bar{\Phi}_2 \bar{\vh},\quad \bar{\vh}
=(\bar{\sph{2}},\cdots,\bar{\sph{n}})^\intercal.
\end{eqnarray}
There exists an operator \(\cal Q\) such that \(\bar{\Phi}_1 \cal
Q\) and \(\bar{\Phi}_2 \cal Q\) are compatible Hamiltonian
operators and \({\cal{ Q}}=T  {\Hop} T^\intercal
{\bar{\Phi}}_{2}^{*}\), where \(\Hop\) is defined in Theorem
\ref{Th1} and \(T\) is defined in the proof of Theorem
\ref{TheoHasimoto}.
\end{theorem}
 \begin{proof} First we notice that \(\bar{\vh}=T \vh\) in Riemannian manifold.
 The Hasimoto transformation is also a Miura transformation
 between \(\vk{t}\) and \(\vkb{t}\), that is,
 \begin{eqnarray*}
 \vk{t}=D_{\vk{}}\vkb{t}=(\Re -\cu )\vh{}
 =D_{\vk{}}(\bar{\Phi}_1 -\cu \bar{\Phi}_2) \bar{\vh}.
 \end{eqnarray*}
 Therefore, \(D_{\vk{}}^{-1}=\bar{\Phi}_2 T\) and \(\Re=D_{\vk{}}\bar{\Phi}_1
 T\). We know that \(\Hop\) and \(\Re \Hop\) are compatible
 Hamiltonian operators. Thus
 \(D_{\vk{}}^{-1} \Hop D_{\vk{}}^{-1,*}=\bar{\Phi}_2 T \Hop T^\intercal
{\bar{\Phi}}_{2}^{*} \) and \(D_{\vk{}}^{-1} \Re \Hop
D_{\vk{}}^{-1,*}=\bar{\Phi}_1 T \Hop T^\intercal
{\bar{\Phi}}_{2}^{*}\)
 are also a compatible Hamiltonian operators. By now, we prove the
 statement.
 \end{proof}
The operator \(\cal Q\) can be computed explicitly and it is
independent on the angles in the Hasimoto transformation.
\renewcommand\ku[2]{{u}_{#1}^{(#2)}}
\renewcommand\kb[2]{\bar{u}_{#1}^{(#2)}}
\section{A Lax pair of vector mKDV}\label{Laxpair}
For the Lie algebra background to find Lax pairs,
see \cite{MR84k:58106,MR1090598}.
The computation of the Lax pair is known in the literature.
In the case of \(\orth{n}\), two derivations, both using different \(\Z/2\)-gradings, are known
to us: \cite{MR2002k:37141,MR1781618}.
Both derivations at some point use some not so obvious steps
and we noted that the steps that were obvious in one were not obvious in the
other. So it seemed like a natural idea to use both \(\Z/2\)-gradings
and do the derivation in a completely obvious way.

We now identify \(\Frame{1}\) with an element in \(\orth{n+1}\)
which we call \(L^{10}\) using the identification
\[
\euc{n}\backepsilon\left(\begin{array}{cccc}0&v^\intercal
\\0& \omega
\end{array}\right)
\stackrel{\phi}{\mapsto} \left(\begin{array}{cccc}0&v^\intercal
\\-v& \omega
\end{array}\right)\in\orth{n+1}
\]
Let \(L=\phi_\star A(\Dx)\), where \(A\) is defined in equation
(\ref{MC}) and let \(\lambda\) be the norm of the tangent vector
to the curve, that is, this vector equals \(\lambda\Frame{1}\).
Then
\begin{eqnarray*}
\lefteqn{L=L^{01}+\lambda L^{10}=}
&&\\
&=&\left(\begin{array}{ccccc}0&0,&0,&\cdots,&0
\\0& 0,&\ku{}{1},&\cdots,& \ku{}{n-1}
\\0& -\ku{}{1},&0,&\cdots,&0
\\\vdots&\vdots&\vdots&\ddots&\vdots
\\0&-\ku{}{n-1},&0,&\cdots,&0
\end{array}\right)
+\lambda\left(\begin{array}{ccccc}0&1,&0,&\cdots,&0
\\-1& 0,&0,&\cdots,&0
\\0& 0,&0,&\cdots,&0
\\\vdots&\vdots&\vdots&\ddots&\vdots
\\0&0,&0,&\cdots,&0
\end{array}\right).
\end{eqnarray*}
In order to get some control over our calculations, we introduce
two \(\Z/2\)-gradings. We do this by partitioning the
\(n+1\)-dimensional matrix into \((1,n)\times(1,n)\)-blocks, and
into \((2,n-1)\times(2,n-1)\)-blocks. We then grade the diagonal
blocks with \(0\) and the off-diagonal blocks with \(1\). This way
we can write \(\orth{n+1}\) as
\[
\orth{n+1}=\lieg^{00}+\lieg^{01}+\lieg^{10}+\lieg^{11}.
\]
We see that \(L^{01}\in\lieg^{01}\) and \(L^{10}\in \lieg^{10}\),
which explains the notation.
Observe that the first grading is consistent with the parity of the power of \(\lambda\),
the second is not. So one does expect the computation to be homogeneous
in the first grading and inhomogeneous with respect to the second.
The second grading codes some useful facts. If \(X\in ker \ ad(L^{10})\),
then \(X\in \lieg^{\cdot 0}\) and if \(X\in im \ ad(L^{10})\),
then \(X\in\lieg^{\cdot 1}\).
If \(X\in \lieg^{10}\), then \(X\in\langle L^{10}\rangle_{\R}\).
Furthermore, \(ad(L^{10})\) gives us a complex structure
on \(\lieg^{01}+\lieg^{10}\).
We now let \(A(\Dt)=M=M_3+\lambda M_2+\lambda^2 M_1 +\lambda^3 M_0\).
The reader may want to experiment with lower order expansions of \(M\).

From the curvature equation (\ref{zero})
\[
F(\Dx,\Dt)=[\Dx+L,\Dt+M]
\]
it follows that the \(\lambda^4\)
term vanishes, that is: \(ad(L^{10})M_0=0\).
We simply choose \(M_0=- L^{10}\).
The next order \(\lambda^3\) gives us
\[
ad(L^{10})M_1=-\Dx M_0-[L^{01},M_0].
\]
or \(M_1+L^{01}\in ker \ ad(L^{10})\).
We take \(M_1=-L^{01}\).

On the \(\lambda^2\) level we find
\[
ad(L^{10})M_2=-\Dx M_1-[L^{01},M_1]=-\Dx L^{01}.
\]
We can solve this by letting  \(M_2=-\Dx ad(L^{10}) L^{01} +\beta L^{10}\),
with \(\beta\) some arbitrary function.
We find on the \(\lambda\) level
\[
ad(L^{10})M_3=\DX{2} ad(L^{10}) L^{01}-\beta_x L^{10}
+[L^{01},\Dx ad(L^{10}) L^{01} ]
-\beta[L^{01},L^{10}]
\]
or
\[
ad(L^{10})M_3=ad(L^{10}) ( \DX{2} L^{01} +\beta L^{01})
-\beta_x L^{10}
-ad(L^{01})ad(\Dx L^{01}) L^{10}
\]
Looking at the grading, we see that
\(ad(L^{01})ad(\Dx L^{01}) L^{10}=\kappa(L^{01},\Dx L^{01}) L^{10}\),
where \(\kappa\) is a bilinear form with values in \(\R\).
We find that on the span of \(L^{10}\) and \(L^{01}\) the operator
\(ad(L^{01})ad(\Dx L^{01})\) behaves like
\[
\langle \vk{},\vk{} \rangle \Dx -\langle \vk{}, \vk{1} \rangle.
\]
In other words, \[
ad(L^{10})M_3=ad(L^{10}) ( \DX{2} L^{01} +\beta L^{01})
-\beta_x L^{10}
+ \langle \vk{}, \vk{1} \rangle L^{10}
\]
So we take \(\beta=\nu+\f2 \langle \vk{}, \vk{} \rangle, \nu\in\R \).
It follows that \(M_3=\DX{2} L^{01} +(\nu+ \f2 \langle \vk{}, \vk{} \rangle)L^{01}
+ X^{00}\).
The last equation (and this gives us the obstruction to solving the
structure equations) reads
\[
F(\Dx,\Dt)=\Dx M_3 -\Dt L^{01} +[L^{01},M_3].
\]
This implies
\begin{eqnarray*}
\Dt L^{01}&=&
\DX{3} L^{01} + \langle \vk{}, \vk{1} \rangle)L^{01}+(\nu+ \f2 \langle \vk{}, \vk{} \rangle)\Dx L^{01}
+\Dx  X^{00}
\\&&
+\Dx[L^{01},\Dx L^{01} ]
+[L^{01}, X^{00}]
-F(\Dx,\Dt)
\end{eqnarray*}
We take \(X^{00}=-[L^{01},\Dx L^{01} ]\)
and obtain
\begin{eqnarray*}
\Dt L^{01}&=&
\DX{3} L^{01} + \langle \vk{}, \vk{1} \rangle L^{01}+(\nu+ \f2 \langle \vk{}, \vk{} \rangle)\Dx L^{01}
\\&&
+ad(L^{01})ad(\Dx L^{01}) L^{01}
-F(\Dx,\Dt)
\\&=&\DX{3} L^{01} + \langle \vk{}, \vk{1} \rangle L^{01}+(\nu+ \f2 \langle \vk{}, \vk{} \rangle)\Dx L^{01}
\\&&
+(\langle \vk{},\vk{} \rangle \Dx -\langle \vk{}, \vk{1} \rangle) L^{01}
-F(\Dx,\Dt)
\\&=&\DX{3} L^{01} +(\nu+ \fr32 \langle \vk{}, \vk{} \rangle)\Dx L^{01}
-F(\Dx,\Dt)
\end{eqnarray*}
Further computation shows that the Killing form of \(L^{01}\)
is given by
\[
K(L^{01},L^{01})=-2(n-2) \langle \vk{}, \vk{} \rangle,
\]
and we can write, at least for \(n>2\),
\begin{eqnarray*}
\Dt L^{01}&=&
\DX{3} L^{01} +(\nu-  \frac{3}{4(n-2)} K(L^{01},L^{01})) \Dx L^{01}
-F(\Dx,\Dt).
\end{eqnarray*}
We now make the choice
\(\Dt=\sum_{i=1}^{n-1} \ku{}{i}
\ku{1}{i}  \Frame{1}+\sum_{i=1}^{n-1}
\ku{1}{i} \Frame{i+1}\),
which is consistent with the condition \(\omega_{11}=0\).
We let \(R_{ijkl}=\Omega_{ij}(\Frame{k},\Frame{l})\).
Constant curvature means that \(R_{ijij}=\cu\) for \(i\neq j\)
and the other \(R\) coefficients are zero.
We obtain
\begin{eqnarray*}
\Dt L^{01}&=&
\DX{3} L^{01} +(\nu-  \frac{3}{4(n-2)} K(L^{01},L^{01})) \Dx L^{01}
-\cu \Dx L^{01}.
\end{eqnarray*}
Observe that the role of the integration constant and the curvature
is the same. This may explain some of the success of the zero-curvature
method.

We can {\em flatten} the equation by absorbing the curvature \(\cu\)
in the connection by taking the integration constant \(\nu\)
equal to \(\cu\):
\begin{eqnarray*}
\Dt L^{01}&=&
\DX{3} L^{01} -  \frac{3}{4(n-2)} K(L^{01},L^{01}) \Dx L^{01}
\end{eqnarray*}
or
\[
\vk{t}=
 \vk{3}+\fr32 ||\vk{}||^2 \vk{1} \mbox{    (vmKDV)}.
\]
From the point of view of integrability this flattening does not
have any influence, since \(\vk{1}\) is a trivial symmetry of the
equation. We remark that while we have found the vmKDV equation
through the Lax pair construction, it also provides a Lax pair for
vmKDV, taking its values in \(\orth{m+1}\), if \(\uu\) is an
\((m-1)\)-vector, cf. \cite{MR2002k:37141}.
\appendix
\section{Geometric structures}\label{OpProof}
In this Appendix we collect the proofs that the operators
\(\mathfrak{I}, \mathfrak{H}\) and \(\Re\) are symplectic,
cosymplectic and hereditary, respectively.
The fact that \(\Re\) is hereditary is equivalent to the compatibility
of the Hamiltonian structures \(\mathfrak{I}^{-1}\) and \(\mathfrak{H}\).

These operators have nonlocal terms by the occurrence of \(\Dxi\).
We refer to \cite{SW00b} for results to on the locality of the generated
symmetries and cosymmetries.
\begin{proposition}\label{SymplecticOp}
The operator \(\mathfrak{I}=\Dx+\uu \Dxi\uu^\intercal\) is symplectic.
\end{proposition}
\begin{proof}
As is proved in \cite{MR94j:58081},
to show that \(\mathfrak{I}\) is symplectic it suffices that
\[
\int ( \langle \mathfrak{I}^\prime [h_1]h_2,h_3\rangle
+\langle \mathfrak{I}^\prime [h_3]h_1,h_2\rangle
+\langle \mathfrak{I}^\prime [h_2]h_3,h_1\rangle )d\pt=0.
\]
One has \(\mathfrak{I}^\prime[h_1]=h_1 \Dxi\uu^\intercal
+ \uu \Dxi h_1^\intercal\). Then
\begin{eqnarray*}
\lefteqn{\langle \mathfrak{I}^\prime [h_1]h_2,h_3\rangle=}&&
\\&=&\langle h_1 \Dxi\uu^\intercal h_2,h_3\rangle
+\langle \uu \Dxi h_1^\intercal h_2,h_3\rangle
\\&=&\Dxi(\langle \uu,h_2\rangle) \langle h_1 ,h_3\rangle
+\Dxi(\langle h_1,h_2\rangle )\langle \uu ,h_3\rangle
\\&\equiv&\Dxi(\langle \uu,h_2\rangle) \langle h_1 ,h_3\rangle
- \Dxi (\langle \uu ,h_3\rangle) \langle h_1,h_2\rangle ,
\end{eqnarray*}
and the result follows upon summation of the \(3\) cyclically permuted terms.
\end{proof}
\begin{proposition}\label{HamiltonianOp}
The operator \(\mathfrak{H}\), given by
\begin{eqnarray*}
\mathfrak{H}=\Dx+\sum_{i=1}^{n-1}\sum_{j=i+1}^n
J_{ij}\uu\Dxi(J_{ij}\uu)^\intercal,
\end{eqnarray*}
where the \(J_{ij}\), given by
\((J_{ij})_{kl}=\delta_{ik}\delta_{jl}-\delta_{il}\delta_{jk}\),
are the generators of \(\orth{n}\),
is a cosymplectic
(or Hamiltonian) operator.
\end{proposition}
\begin{proof}
We prove this by checking the conditions of
    Theorem 7.8 and Corollary 7.21 in \cite{MR94g:58260}.

The associated bi-vector of \(\Hop\) is by definition
\begin{eqnarray*}
\Theta_{\Hop}
&=&\f2 \int \left(  \theta \dot{\wedge}\Hop \theta \right) d\pt\\
\\&=&\f2 \int \left(  \theta \dot{\wedge}\left(\theta_1+\sum_{i=1}^{n-1}\sum_{j=i+1}^n J_{ij}\uu\Dxi(J_{ij}\uu)^\intercal\theta\right)\right) d\pt\\
\\&=&\f2 \int \theta \dot{\wedge}\theta_1 d\pt
+\f2 \sum_{i=1}^{n-1}\sum_{j=i+1}^n\int  \langle J_{ij}\uu,\theta \rangle\wedge\langle J_{ij}\uu,\theta\rangle_{-1} d\pt,
\end{eqnarray*}
where \(\theta=(\vartheta^1, \cdots,\vartheta^n)\) and
\(\vartheta_i^j=\frac{\partial^i \vartheta^j}{\partial x^i}\), etc.
Here \(\dot{\wedge}\) means that one needs to take the
ordinary inner product between
the vectors \(\theta \) and \( \Hop \theta \).
The elements of these vectors are then multiplied using the
ordinary wedge product.
We need to check the vanishing of the Schouten bracket \([\Hop,\Hop]\)
which is equivalent to the Jacobi identity for the Lie bracket defined
by \(\Hop\).
\begin{eqnarray*}
\lefteqn{[\Hop,\Hop]=\Pr \V_{\Hop \theta}(\Theta_{\Hop})=}&&\\
\\&=&\sum_{i<j} \int \left(
\langle J_{ij}\Hop \theta,\theta \rangle\wedge
\langle J_{ij}\uu,\theta\rangle_{-1}  \right) d\pt \\
\\&=&\sum_{i<j} \int
 J_{ij}\theta_1\dot{\wedge}\theta \wedge
\langle J_{ij}\uu,\theta\rangle_{-1}  d\pt \\
\\&+&\sum_{i<j} \sum_{k<l}\int
\langle J_{kl}\uu,\theta\rangle_{-1}
J_{ij}J_{kl}\uu
\dot{\wedge} \theta \wedge
\langle J_{ij}\uu,\theta\rangle_{-1}  d\pt.
\end{eqnarray*}
The first term is zero because
\begin{eqnarray*}
\lefteqn{J_{ij}\theta_1\dot{\wedge}\theta \wedge
\langle J_{ij}\uu,\theta\rangle_{-1} =}&
\\&=&(\theta_1^j\wedge\theta^i-\theta_1^i\wedge\theta^j) \wedge
\Dxi (\uu^j\theta^i -\uu^i\theta^j)
\\&=&(\theta_1^j\wedge\theta^i+\theta^j\wedge\theta_1^i) \wedge
\Dxi (\uu^j\theta^i -\uu^i\theta^j)
\\&\equiv&
-\theta^j\wedge\theta^i \wedge
(\uu^j\theta^i -\uu^i\theta^j)
\\&=&0.
\end{eqnarray*}
The second term is zero because
\begin{eqnarray*}
\lefteqn{
\sum_{i<j} \sum_{k<l}
\langle J_{kl}\uu,\theta\rangle_{-1}
\wedge \langle J_{ij}J_{kl}\uu
,\theta\rangle  \wedge
\langle J_{ij}\uu,\theta\rangle_{-1}  }&
\\&=&
\f2\sum_{i<j} \sum_{k<l}
\langle J_{kl}\uu,\theta\rangle_{-1}
\wedge \langle[ J_{ij},J_{kl}]\uu
,\theta\rangle  \wedge
\langle J_{ij}\uu,\theta\rangle_{-1}
\\&=&
\sum_{i<j,k<l,i<k}
\langle J_{kl}\uu,\theta\rangle_{-1}
\wedge \langle\delta_{lj}J_{ik}\uu
,\theta\rangle  \wedge
\langle J_{ij}\uu,\theta\rangle_{-1}
\\&-&\sum_{i<j,k<l,i<l}
\langle J_{kl}\uu,\theta\rangle_{-1}
\wedge \langle\delta_{jk}J_{il}\uu
,\theta\rangle  \wedge
\langle J_{ij}\uu,\theta\rangle_{-1}
\\&+&\sum_{i<j,k<l,j<l}
\langle J_{kl}\uu,\theta\rangle_{-1}
\wedge \langle\delta_{ik}J_{jl}\uu
,\theta\rangle  \wedge
\langle J_{ij}\uu,\theta\rangle_{-1}
\\&=&
\sum_{i<j<k} \Dx\left(
\langle J_{jk}\uu,\theta\rangle_{-1}
\wedge \langle J_{ij}\uu
,\theta\rangle_{-1}  \wedge
\langle J_{ik}\uu,\theta\rangle_{-1} \right)
\\&\equiv&0.
\end{eqnarray*}
Thus the result follows.
\end{proof}
\begin{proposition}\label{HereditaryOp}
The operator
\begin{eqnarray*}
\Re&=&\Dx^2+\langle \vk{},\ \vk{}\rangle +\vk{1} \Dxi \langle \vk{},\ \cdot\rangle
-\sum_{i<j}^{n} J_{ij} \vk{} \Dxi
\langle J_{ij} \vk{1},\ \cdot\rangle
\end{eqnarray*}
is hereditary.
\end{proposition}
\begin{proof}
Notice that for the matrices \(J_{ij}\) we have
\begin{eqnarray}\label{JR}
\sum_{i<j}^{n} J_{ij} \vk{} \Dxi
\langle J_{ij} \vk{1},\ P\rangle
= (\Dxi P \vk{1}^\intercal)\vk{}  - (\Dxi \vk{1} P^\intercal) \vk{}.
\end{eqnarray}
We recall that an hereditary operator \(\Re\) is characterized by
the property that
\begin{eqnarray*}
\Re D_{\Re}[P] (Q)-D_{\Re}[\Re P](Q) \quad \mbox{is symmetric with
respect to vectors \(P\) and \(Q\), cf. \cite{MR84j:58046}.}
\end{eqnarray*}
We first compute \(D_{\Re}[P] (Q)\), meanwhile we drop the terms that are
symmetric with respect to \(P\) and \(Q\) by using "\(\equiv\)".
\begin{eqnarray*}
D_{\Re}[P] (Q)&=&
2 \langle \vk{},\ P\rangle Q +P_x \Dxi \langle \vk{},\ Q\rangle +\vk{1} \Dxi \langle P,\ Q\rangle\\
&&-\sum_{i<j}^{n} J_{ij} P \Dxi
\langle J_{ij} \vk{1},\ Q\rangle
-\sum_{i<j}^{n} J_{ij} \vk{} \Dxi
\langle J_{ij} P_x,\ Q\rangle
\\
&\equiv&
\langle \vk{}, P\rangle Q +P_x \Dxi \langle \vk{},\ Q\rangle
-\sum_{i<j}^{n} J_{ij} P \Dxi
\langle J_{ij} \vk{1},\ Q\rangle,
\end{eqnarray*}
where we use the relation \(\sum_{i<j}^{n} J_{ij} \vk{} \Dxi
\langle J_{ij} P_x,\ Q\rangle
\equiv \langle \vk{},\ P\rangle Q\), which can be proved by applying formula (\ref{JR}).
Notice that
\[\Re P=P_{xx}+ \langle \vk{},\ \vk{}\rangle P
+\vk{1} \Dxi \langle \vk{},\ P\rangle
-\sum_{i<j}^{n} J_{ij} \vk{} \Dxi
\langle J_{ij} \vk{1},\ P\rangle
\]
and the terms in \(\Re D_{\Re}[P] (Q)-D_{\Re}[\Re P](Q)\) are \(\vk{}\)--degree of \(1\)
or \(3\). Its first degree terms are
\begin{eqnarray*}
&&\langle \vk{2},\ P\rangle  Q
+2 \langle \vk{1},\ P_x\rangle  Q +\langle \vk{},\ P_{xx}\rangle  Q
\\&&
+2 \langle \vk{1},\ P\rangle  Q_x +2 \langle \vk{},\ P_x\rangle  Q_x +\langle \vk{},\ P\rangle  Q_{xx}
\\&&
+P_x \langle \vk{},\ Q_x\rangle  +P_x \langle \vk{1},\ Q\rangle  +2 P_{xx} \langle \vk{},\ Q\rangle
+P_{xxx} \Dxi \langle \vk{},\ Q\rangle
\\&&
-\sum_{i<j}^{n} J_{ij} P
\langle J_{ij} \vk{2},\ Q\rangle
-\sum_{i<j}^{n} J_{ij} P
\langle J_{ij} \vk{1},\ Q_x\rangle
\\&&
-2\sum_{i<j}^{n} J_{ij}P_x
\langle J_{ij} \vk{1},\ Q\rangle
-\sum_{i<j}^{n} J_{ij} P_{xx} \Dxi
\langle J_{ij} \vk{1},\ Q\rangle
\\&&
-2 \langle P_{xx},\ \vk{}\rangle Q - P_{xxx} \Dxi \langle \vk{},\ Q\rangle  -\vk{1} \Dxi \langle P_{xx},\ Q\rangle
\\&&
+\sum_{i<j}^{n} J_{ij}  P_{xx} \Dxi
\langle J_{ij} \vk{1},\ Q\rangle
+\sum_{i<j}^{n} J_{ij} \vk{} \Dxi
\langle J_{ij} P_{xxx},\ Q\rangle
\\
&\equiv&0.
\end{eqnarray*}
Similarly, we can check its third degree terms are also equivalent to zero.
Therefore,
\begin{eqnarray*}
\left(\Re D_{\Re}[P] (Q)\right)^j
-\left(D_{\Re}[\Re P](Q)\right)^j \equiv 0,
\end{eqnarray*}
that is, the operator \(\Re\) is hereditary.
\end{proof}
\section{Generalized Hasimoto transformation}\label{Hasimoto}
\renewcommand\vk[1]{{\bf \bar{u}}_{#1}}
\renewcommand\kb[2]{u_{#1}^{(#2)}}
\renewcommand\ku[2]{{\bar{u}}_{#1}^{(#2)}}
\renewcommand\bFrame[1]{\mathbf{e}_{#1}}
\renewcommand\Frame[1]{\bar{\mathbf{e}}_{#1}}
In this Appendix we complete the proof of Theorem
\ref{TheoHasimoto}. The formula contain possibly non-existing
angles. If this is the case the convention is to take the first
term in the product to be equal to one.
\begin{lemma}
The matrix \(T\) gauges \(\suboo\) into zero, cf. the proof of
Theorem \ref{TheoHasimoto} for the notations of \(T\) and
\(\suboo\).,that is \(\suboo T=\Dx T\).
\end{lemma}
\begin{proof}

Notice that \(R=R_{n-1,n}\cdots R_{2n}R_{n-2,n-1}\cdots R_{34}R_{24}R_{23}\) since
\(R_{ij}\) commutes \(R_{kl}\) when \(\{i,j\}\bigcap \{k,l\}=\emptyset\).
Let \(R^{(i)}=R_{i-1,i}\cdots R_{3i} R_{2i}\) and the part without the first row
and the first column is denoted by \(\BB{i}{}\). Then one writes
\(R=R^{(n)}R^{(n-1)}\cdots R^{(4)} R^{(3)}\)
and \(T=\BB{n}{}\BB{n-1}{}\cdots \BB{4}{} \BB{3}{}\).
The matrix \(\BB{i}{}\) is generated by
\begin{eqnarray*}
\left\{\begin{array}{lll}
\BB{i}{jj}=1,\quad j>i-1\\
\BB{i}{jj}= \cos \theta_{j+1,i}, \quad 1<j<i-1,\\
\BB{i}{1,i-1}=\sin\theta_{2i}\\
\BB{i}{lj}=\sin \theta_{l+1,i} \frac{\partial
\BB{i}{l-1,j}}{\partial \theta_{l,i}},\quad j<l<i-1\\
\BB{i}{i-1,j}=\frac{\partial
\BB{i}{i-2,j}}{\partial \theta_{i-1,i}}.
\end{array}\right.
\end{eqnarray*}
and the rest entries are equal to zero.

We first perform the gauge transformation of \(\BB{n}{}\)  on
\(\suboo\). Notice the last column of \(\suboo \BB{n}{}-\Dx
\BB{n}{}\) has to be zero since the others do not affect it. This
implies that
\begin{eqnarray}
&&\ku{}{2} \BB{n}{2,n-1}-\Dx \BB{n}{1,n-1}=0;\label{AFormu2}\\
&&\ku{}{k}\BB{n}{k,n-1}-
\ku{}{k-1} \BB{n}{k-2,n-1}-\Dx \BB{n}{k-1,n-1}=0,
\quad 2<k\leq n-1;\label{AFormuk}\\
&&\ku{}{n-1} \BB{n}{n-2,n-1}+\Dx \BB{n}{n-1,n-1}=0,\label{AFormun}
\end{eqnarray}
where \(\BB{n}{j,n-1}=\sin\theta_{j+1,n} \cos\theta_{j,n}\cdots \cos\theta_{2n}\)
and \(1\leq j\leq n-1\).
The formula (\ref{AFormu2}) immediately leads to
\begin{eqnarray*}
\ku{}{2}=\frac{\Dx \theta_{2n}}{\sin \theta_{3n}}.
\end{eqnarray*}
Now we prove the following formula by induction on \(i\).
\begin{eqnarray}\label{AFormui}
&&\ku{}{i}=\frac{\Dx \theta_{i,n}}{\sin \theta_{i+1,n}}
+\frac{\cos \theta_{i,n}}{\sin \theta_{i+1,n}}
\frac{\sin \theta_{i-1,n}}{\cos \theta_{i-1,n}} \ku{}{i-1},\quad 2<i\leq n-1.
\end{eqnarray}
Using (\ref{AFormuk}) and assuming that (\ref{AFormui})
is true for \(2<l\leq k-1\), we compute
\begin{eqnarray*}
&&\ku{}{k} \BB{n}{k,n-1}=
\ku{}{k} \sin\theta_{k+1,n}\cos\theta_{k,n} \cdots \cos\theta_{2,n}\\
&=&\ku{}{k-1} \BB{n}{k-2,n-1}+\Dx \BB{n}{k-1,n-1}\\
&=&\ku{}{k-1}\sin\theta_{k-1,n}\cos\theta_{k-2,n} \cdots \cos\theta_{2,n}\\
&&+ \cos\theta_{k,n}\cos\theta_{k-1,n} \cdots \cos\theta_{2,n} \Dx \theta_{k,n}\\
&&-\sum_{l=2}^{k-1} \sin\theta_{k,n}\cos\theta_{k-1,n} \cdots
\sin\theta_{l,n}\cdots \cos\theta_{2,n}\Dx \theta_{l,n}\\
&=&\ku{}{k-1}\cos^2 \theta_{k,n} \sin\theta_{k-1,n}\cos\theta_{k-2,n}
\cdots \cos\theta_{2,n}\\
&&+ \cos\theta_{k,n}\cos\theta_{k-1,n} \cdots \cos\theta_{2,n} \Dx \theta_{k,n}.
\end{eqnarray*}
Thus, we proved the formula is valid for \(\ku{}{k}\).
By directly computation, we can check the identity (\ref{AFormun}).

Next, we prove \(\suboo \BB{n}{}-\Dx \BB{n}{}=\BB{n}{} C\), where
\(C\) is an anti-symmetric matrix with
\begin{eqnarray*}
C_{k-1,k}=-C_{k,k-1}= \sum_{l=2}^{k}
\frac{\cos \theta_{k+1,n}}{\sin \theta_{k+1,n}}
\frac{1}{\sin \theta_{k,n}}
\frac{\sin \theta_{l,n}}{\cos \theta_{l,n}} \Dx\theta_{l,n},\quad 2\leq i\leq n-2,
\end{eqnarray*}
and the remaining entries are equal to zero.

To do so, we first rewrite (\ref{AFormui}) as
\begin{eqnarray}\label{AFormud}
&&\ku{}{i}= \sum_{l=2}^{i} \frac{1}{\sin \theta_{i+1,n}}
\frac{\cos \theta_{i,n}}{\sin \theta_{i,n}}
\frac{\sin \theta_{l,n}}{\cos \theta_{l,n}} \Dx\theta_{l,n},
\quad 2\leq i\leq n-1.
\end{eqnarray}
Let us check that
\begin{eqnarray}\label{FFm}
0&=&\ku{}{i+1}\BB{n}{i+1,j}-
\ku{}{i} \BB{n}{i-1,j}-\Dx\BB{n}{i,j}\nonumber\\
&-&\BB{n}{i,j-1} C_{j-1,j}+\BB{n}{i,j+1} C_{j,j+1}.
\end{eqnarray}
It is easy to see that this formula is valid for \(j>i+1\) since
\(\BB{n}{l,k}=0\) when \(l<k\) and \(k\neq n-1\).
When \(j=i+1\), we have
\begin{eqnarray*}
\ku{}{i+1}\BB{n}{i+1,i+1} -\BB{n}{i,i} C_{i,i+1}
=\ku{}{i+1}\cos \theta_{i+2,n} -\cos \theta_{i+1,n} C_{i,i+1}=0.
\end{eqnarray*}
Similarly, for \(j=i\), we can check
\(\ku{}{i+1}\BB{n}{i+1,i}-\Dx\BB{n}{i,i} -\BB{n}{i,i-1} C_{i-1,i}=0.\)
Now we concentrate on the case \(j<i\). Notice that
\begin{eqnarray*}
\BB{n}{i,j}=-\sin \theta_{i+1,n} \cos \theta_{i,n} \cdots
\cos \theta_{j+2,n} \sin \theta_{j+1,n}, \quad i>j.
\end{eqnarray*}
So, we have
\begin{eqnarray*}
&&\ku{}{i+1}\BB{n}{i+1,j}-
\ku{}{i} \BB{n}{i-1,j}-\Dx\BB{n}{i,j}
-\BB{n}{i,j-1} C_{j-1,j}+\BB{n}{i,j+1} C_{j,j+1}\\
&&=\sum_{l=2}^{i+1} \frac{1}{\sin \theta_{i+2,n}}
\frac{\cos \theta_{i+1,n}}{\sin \theta_{i+1,n}}
\frac{\sin \theta_{l,n}}{\cos \theta_{l,n}} \Dx\theta_{l,n}
(-\sin \theta_{i+2,n} \prod_{k=j+2}^{i+1}\cos \theta_{k,n}
\sin \theta_{j+1,n})
\\
&&-\sum_{l=2}^{i} \frac{1}{\sin \theta_{i+1,n}}
\frac{\cos \theta_{i,n}}{\sin \theta_{i,n}}
\frac{\sin \theta_{l,n}}{\cos \theta_{l,n}} \Dx\theta_{l,n}
(-\sin \theta_{i,n} \prod_{k=j+2}^{i-1}\cos \theta_{k,n}
\sin \theta_{j+1,n})
\\
&&-\Dx \left(-\sin \theta_{i+1,n}\prod_{k=j+2}^{i}\cos \theta_{k,n}
\sin \theta_{j+1,n}\right)
\\
&&-(-\sin \theta_{i+1,n} \prod_{k=j+1}^{i}\cos \theta_{k,n}
\sin \theta_{j,n})
\sum_{l=2}^{j}
\frac{\cos \theta_{j+1,n}}{\sin \theta_{j+1,n}}
\frac{1}{\sin \theta_{j,n}}
\frac{\sin \theta_{l,n}}{\cos \theta_{l,n}} \Dx\theta_{l,n}
\\
&&+(-\sin \theta_{i+1,n} \prod_{k=j+3}^{i}\cos \theta_{k,n}
\sin \theta_{j+2,n})
\sum_{l=2}^{j+1}
\frac{\cos \theta_{j+2,n}}{\sin \theta_{j+2,n}}
\frac{1}{\sin \theta_{j+1,n}}
\frac{\sin \theta_{l,n}}{\cos \theta_{l,n}} \Dx\theta_{l,n}
\\
&&=0.
\end{eqnarray*}

Now we apply induction procedure on the matrix \(C\) since it is of
the form of Fren{\^e}t frame.  Using (\ref{AFormud}) for the \(n-1\) case, we have
\begin{eqnarray*}
&&C_{k-1,k}=\sum_{l=2}^{k} \frac{1}{\sin \theta_{k+1,n-1}}
\frac{\cos \theta_{k,n-1}}{\sin \theta_{k,-1n}}
\frac{\sin \theta_{l,n-1}}{\cos \theta_{l,n-1}} \Dx\theta_{l,n-1},
\quad 2\leq k\leq n-2.
\end{eqnarray*}
Form these, we can solve \(\Dx \theta_{k,n}\) for \(2\leq k\leq n-2\).
\begin{eqnarray*}
\Dx\theta_{k,n}&=&\sum_{l=2}^{k}
\frac{\sin \theta_{k+1,n}}{\cos \theta_{k+1,n}}
\frac{\cos \theta_{k,n}}{\sin \theta_{k+1,n-1}}
\frac{\cos \theta_{k,n-1}}{\sin \theta_{k,n-1}}
\frac{\sin \theta_{l,n-1}}{\cos \theta_{l,n-1}} \Dx\theta_{l,n-1}\\
&&-
\sum_{l=2}^{k-1}
\frac{\cos \theta_{k,n}}{\sin \theta_{k,n}}
\frac{\sin \theta_{l,n}}{\cos \theta_{l,n}} \Dx\theta_{l,n}
\quad 2\leq k\leq n-2.
\end{eqnarray*}
Assume that the formula (\ref{AFormth}) is valid for \(j= n-1\),
and \(2\leq i\leq n-3\). When \(k=2\), we obtain
\begin{eqnarray*}
\Dx\theta_{2,n}&=&
\frac{\sin \theta_{3,n}}{\cos \theta_{3,n}}
\frac{\cos \theta_{2,n}}{\sin \theta_{3,n-1}}
\prod_{l=4}^{n-1}\frac{\cos \theta_{2,l}}{\cos \theta_{3,l}}
\sin \theta_{3,n-1} \subo_2\\
&=&
\prod_{l=4}^{n}\frac{\cos \theta_{2,l}}{\cos \theta_{3,l}}
\sin \theta_{3,n} \subo_2.
\end{eqnarray*}
Now assuming that the formula (\ref{AFormth}) is valid for \(j= n\),
and \(2\leq i\leq k-1\), we compute
\begin{eqnarray*}
\Dx\theta_{k,n}&=&
\sum_{l=2}^{k}
\frac{\sin \theta_{k+1,n}}{\cos \theta_{k+1,n}}
\frac{\cos \theta_{k,n}}{\sin \theta_{k+1,n-1}}
\frac{\cos \theta_{k,n-1}}{\sin \theta_{k,n-1}}
\frac{\sin \theta_{l,n-1}}{\cos \theta_{l,n-1}}\cdot \\
&&\left(
\prod_{r=l+2}^{n-1}\frac{\cos \theta_{lr}}{\cos \theta_{l+1,r}}
\sin \theta_{l+1,n-1} \subo_l
-\prod_{r=l+1}^{n-2}\frac{\cos \theta_{l-1,r}}{\cos \theta_{l,r}}
\sin \theta_{l-1,n-1} \subo_{l-1}
\right)
\\
&-&
\sum_{l=2}^{k-1}
\frac{\cos \theta_{k,n}}{\sin \theta_{k,n}}
\frac{\sin \theta_{l,n}}{\cos \theta_{l,n}} \cdot\\
&& \left(
\prod_{r=l+2}^n\frac{\cos \theta_{lr}}{\cos \theta_{l+1,r}}
\sin \theta_{l+1,n} \subo_l
-\prod_{r=l+1}^{n-1}\frac{\cos \theta_{l-1,r}}{\cos \theta_{l,r}}
\sin \theta_{l-1,n} \subo_{l-1},
\right)\\
&=&
\prod_{r=k+2}^{n}\frac{\cos \theta_{kr}}{\cos \theta_{k+1,r}}
\sin \theta_{k+1,n} \subo_k -
\prod_{r=k+1}^{n-1}\frac{\cos \theta_{k-1,r}}{\cos \theta_{k,r}}
\sin \theta_{k-1,n} \subo_{k-1}.
\end{eqnarray*}
Therefore, (\ref{AFormth}) is also valid for \(j=n\) and \(i=k\).

Substituting (\ref{AFormth}) into (\ref{AFormud}), we have
\begin{eqnarray*}
\ku{}{i}&=& \sum_{l=2}^{i} \frac{1}{\sin \theta_{i+1,n}}
\frac{\cos \theta_{i,n}}{\sin \theta_{i,n}}
\frac{\sin \theta_{l,n}}{\cos \theta_{l,n}} \cdot\\
&&
\left(
\prod_{r=l+2}^{n}\frac{\cos \theta_{lr}}{\cos \theta_{l+1,r}}
\sin \theta_{l+1,n} \subo_l -
\prod_{r=l+1}^{n-1}\frac{\cos \theta_{l-1,r}}{\cos \theta_{l,r}}
\sin \theta_{l-1,n} \subo_{l-1}
\right)\\
&=&
\prod_{r=i+2}^{n}\frac{\cos \theta_{ir}}{\cos \theta_{i+1,r}} \subo_i.
\end{eqnarray*}
Finally, we prove (\ref{AFormA}) for \(i=n-1\).
Taking \(i=n-1\) in (\ref{AFormui}), we have
\begin{eqnarray*}
&&\ku{}{n-1}=\Dx \theta_{n-1,n}
+\cos \theta_{n-1,n}
\frac{\sin \theta_{n-2,n}}{\cos \theta_{n-2,n}} \ku{}{n-2}
\end{eqnarray*}
On the other hand, we have \(\ku{}{n-1}=\subo_{n-1}\)
and \(\ku{}{n-2}= \frac{\cos \theta_{n-2,n}}{\cos \theta_{n-1,n}} \subo_{n-2}.\)
Therefore
\begin{eqnarray*}
\subo_{n-1}=\Dx \theta_{n-1,n} + \sin \theta_{n-2,n} \subo_{n-2}.
\end{eqnarray*}
By now we completely proved the statement.
\end{proof}

\end{document}